\documentclass[11pt, a4paper]{article}
\usepackage{amsmath}
\usepackage{wasysym}
\usepackage{latexsym}
\usepackage{amsfonts}
\usepackage{mathrsfs}
\usepackage{amssymb}
\usepackage{ifsym}
\usepackage{dsfont}
\usepackage[all]{xy}
\usepackage{authblk}

\newtheorem{Definitions1}{Definition}[section]

\newtheorem{Theorems1}{Theorem}[section]

\newtheorem{Coroll1}[Theorems1]{Corollary}
\newtheorem{Lemma1}[Theorems1]{Lemma}

\newtheorem{Examp1}{Example}[section]

\newtheorem{Quest1}{Question}[section]

\newenvironment{proof}[1][Proof]{\begin{trivlist}
\item[\hskip \labelsep {\bfseries #1}]}{\end{trivlist}}

\newcommand{\qed}{\nobreak \ifvmode \relax \else
      \ifdim\lastskip<1.5em \hskip-\lastskip
      \hskip1.5em plus0em minus0.5em \fi \nobreak
      \vrule height0.75em width0.5em depth0.25em\fi}
      
\newcommand{\dotminus}{
\begin{picture}(12,6)(-4,-2) 
  \put(-1, 0){ \makebox(0,0){$-$} }
  \put(-1, 2){ \makebox(0,0){$\cdot$} } 
  \end{picture}}

\begin{document}
\title{The subset relation and $2$-stratified sentences in set theory and class theory}
\author{Zachiri McKenzie}
\affil{\texttt{zach.mckenzie@gmail.com}}
\maketitle

\begin{abstract}
Hamkins and Kikuchi (2016 and 2017) show that in both set theory and class theory the definable subset ordering of the universe interprets a complete and decidable theory. This paper identifies the minimum subsystem of $\mathrm{ZF}$, $\mathrm{BAS}$, that ensures that the definable subset ordering of the universe interprets a complete theory, and classifies the structures that can be realised as the subset relation in a model of this set theory. Extending and refining Hamkins and Kikuchi's result for class theory, a complete extension, $\mathrm{IABA}_{\mathrm{Ideal}}$, of the theory of infinite atomic boolean algebras and a minimum subsystem, $\mathrm{BAC}^+$, of $\mathrm{NBG}$ are identified with the property that if $\mathcal{M}$ is a model of $\mathrm{BAC}^+$, then $\langle M, \mathcal{S}^\mathcal{M}, \subseteq^\mathcal{M} \rangle$ is a model of $\mathrm{IABA}_{\mathrm{Ideal}}$, where $M$ is the underlying set of $\mathcal{M}$, $\mathcal{S}^\mathcal{M}$ is the unary predicate that distinguishes sets from classes and $\subseteq^\mathcal{M}$ is the definable subset relation. These results are used to show that that $\mathrm{BAS}$ decides every $2$-stratified sentence of set theory and $\mathrm{BAC}^+$ decides every $2$-stratified sentence of class theory.   
\end{abstract}

\section[Introduction]{Introduction}

Recent work by Hamkins and Kikuchi \cite{hk16, hkxx} initiates the study of the structure $\langle M, \subseteq^\mathcal{M} \rangle$ where $\subseteq^{\mathcal{M}}$ is the definable subset relation in a model of set theory or class theory, $\mathcal{M}$, with domain $M$. \cite{hk16} shows that in models of set theory the definable subset relation is an atomic unbounded relatively complemented distributive lattice, which is a complete and decidable theory. In \cite{hkxx}, it is shown that if $\mathcal{M}=\langle M, \in^\mathcal{M} \rangle$ is a model of set theory that satisfies some very mild set-theoretic conditions (conditions that are satisfied in all models of $\mathrm{ZFC}$ and even in all nonstandard models of finite set theory), then the structure $\langle M, \subseteq^\mathcal{M} \rangle$ is $\omega$-saturated. Thus, if $\mathcal{M}$ is countable, then the structure $\langle M, \subseteq^\mathcal{M} \rangle$ is unique up to isomorphism. In addition, \cite{hkxx} also shows that, in any model of von Neumann-Bernays-G\"{o}del class theory plus the class formulation of the axiom choice the definable subset relation is an infinite atomic Boolean algebra that is $\omega$-saturated. Thereby showing that, as is the case with set theory, the theory of the subset relation in class theory is a complete decidable theory whose structure, for countable models of class theory, is unique. These results are presented as evidence that the subset relation alone is insufficient to serve as a foundation for mathematics. Indeed, the fact that the theory of the definable subset relation in set theory and class theory is a complete theory shows that, after some point, varying the axioms of set theory or class theory in ways that we know fundamentally change the mathematics that can be done (for example, adding or removing the Axiom of Infinity or the Axiom of Choice, or varying the amount of comprehension or collection that is available) does not alter the first order theory of the subset ordering of the universe. Similarly, it also shows that whatever fragment of set theory that is expressible in the language that only consists of the subset relation must be too weak to be able to express any assertion that is not decided by the weak fragment of set theory or class theory that fixes the complete theory of the definable subset ordering. This paper makes these observations precise by identifying the minimal subsystems of set theory and class theory that fix the complete theories of the subset ordering identified in \cite{hk16, hkxx}. We also identify the precise fragment of set theory and class theory that is expressible using only the subset relation, thereby showing that weak subsystems of set theory and class theory decide every sentence in this fragment.

The weak system {\it Adjunctive Set Theory with Boolean operations} ($\mathrm{BAS}$) is axiomatised by extensionality, emptyset, axioms asserting that for all sets $x$ and $y$, the sets $x \cap y$, $x \cup y$, $x-y$ and $\{x\}$ exist, and an axiom asserting that there is no universal set. In section \ref{Sec:SetTheoryAndMer}, it is shown that if $\mathcal{M}= \langle M, \in^\mathcal{M} \rangle$ is a model of $\mathrm{BAS}$ then $\langle M, \subseteq^\mathcal{M} \rangle$ is an atomic unbounded relatively complemented distributive lattice with the same number of elements as atoms. Conversely, we show that every atomic unbounded relatively complemented distributive lattice with the same number of atoms as elements can be realised as the definable subset relation of a model of $\mathrm{BAS}$. It follows that if $\mathcal{M}= \langle M, \in^\mathcal{M} \rangle$ is a model of Tarski's Adjunctive Set Theory, then $\langle M, \subseteq^\mathcal{M} \rangle$ is an atomic unbounded relatively complemented distributive lattice if and only if $\mathcal{M}$ satisfies $\mathrm{BAS}$.

In section \ref{Sec:2StratifiedSentences}, we identify the $2$-stratified sentences as the exact fragment of the language of set theory that is expressible in the structure $\langle M, \subseteq^\mathcal{M} \rangle$, where $\mathcal{M}$ is a model of set theory with domain $M$ and $\subseteq^\mathcal{M}$ is the definable subset relation. We present a translation, $\tau$, of $2$-stratified formulae of set theory into formulae in language of orderings such that if $\mathcal{M}$ is a model of set theory with domain $M$ and $\phi$ is a $2$-stratified sentence in the language of set theory, then $\mathcal{M} \models \phi$ if and only if $\langle M, \subseteq^\mathcal{M} \rangle \models \phi^\tau$. Since the theory of the definable subset relation in the theory $\mathrm{BAS}$ is complete, this shows that any extension of $\mathrm{BAS}$ decides every $2$-stratified sentence in the language of set theory. We show that this result is optimal by expressing a version of the Axiom of Choice using a $3$-stratified sentence of set theory.

Section \ref{Sec:ClassTheoryAndIABA} turns to investigating the definable subset relation in models of class theory. The system {\it Adjunctive Class Theory with Boolean operations} ($\mathrm{BAC}$) consists of extensionality for classes, axioms asserting that the sets are closed downwards under the subset relation and satisfy the theory $\mathrm{BAS}$, and axioms asserting that for all classes $X$ and $Y$, the classes $X \cap Y$, $X \cup Y$ and the complement of $X$ exist. We show that if $\mathcal{M}$ is a model of $\mathrm{BAC}$ with domain $M$ and sets $\mathcal{S}^\mathcal{M}$, then $\langle M, \subseteq^\mathcal{M} \rangle$ is an infinite atomic Boolean algebra and $\mathcal{S}^\mathcal{M}$ is a proper ideal of $\langle M, \subseteq^\mathcal{M} \rangle$ that is the same size as, and contains all of, the atoms of $\langle M, \subseteq^\mathcal{M} \rangle$. Conversely, every infinite atomic Boolean algebra can be realised as the subset ordering of a model $\mathcal{M}$ of $\mathrm{BAC}$. While the translation, $\tau$, introduced in section \ref{Sec:2StratifiedSentences} shows that $\mathrm{BAC}$ decides every $2$-stratified sentence in the language of set theory (the language that only includes $\in$), we note that there is a $2$-stratified sentence of the language of class theory including a unary predicate distinguishing the sets that is independent of $\mathrm{BAC}$. Using a complete decidable theory studied in \cite{dm13}, which in this paper we call $\mathrm{IABA}_{\mathrm{Ideal}}$, we identify an extension $\mathrm{BAC}^+$ of $\mathrm{BAC}$ such that for models $\mathcal{M}$ of this theory, the structure $\langle M, \mathcal{S}^\mathcal{M}, \subseteq^\mathcal{M} \rangle$ satisfies $\mathrm{IABA}_{\mathrm{Ideal}}$, where $M$ is the underlying set of $\mathcal{M}$, $\subseteq^\mathcal{M}$ is the definable subset relation and $\mathcal{S}^\mathcal{M}$ distinguishes the sets of $\mathcal{M}$. This shows that the theory $\mathrm{BAC}^+$, a subsystem of von Neumann-Bernays-G\"{o}del class theory, decides every $2$-stratified sentence in the language of class theory including a unary predicate that distinguishes the sets.   

The results in this paper are influenced by the work of Gri\u{s}hin \cite{gri73} on subsystems of Quine's `New Foundations' Set Theory ($\mathrm{NF}$).  Gri\u{s}hin shows that the Simple Theory of Types restricted to only two types ($\mathrm{TST}_2$) is a complete and decidable theory with essentially the same expressive power as the theory of infinite atomic Boolean algebras. The results of section \ref{Sec:2StratifiedSentences} of this paper are an analogue of this result of Gri\u{s}hin's for set theories that refute the existence of a universal set. Gri\u{s}hin also shows that the fragment of Quine's `New Foundations' Set Theory axiomatised by extensionality and all $2$-stratifiable instances of the comprehension scheme ($\mathrm{NF}_2$) is finitely axiomatised by extensionality and axioms asserting that there exists a universal set ($V$) and for all sets $x$ and $y$, $x \cap y$, $x \cup y$, $\{x\}$ and $V-x$. Moreover, models of $\mathrm{NF}_2$ can be obtained from any infinite atomic Boolean algebra with the same number of atoms as elements. The set theory $\mathrm{BAS}$ introduced in this paper is the theory $\mathrm{NF}_2$ with the axiom asserting the existence of a universal set replaced with its negation.\medskip

\noindent {\bf Acknowledgements:} I am grateful to Randall Holmes for a discussion that led to these results, and Ruizhi Yang, Thomas Forster and an anonymous referee for helpful comments on earlier versions of this paper.

\section[Background]{Background} \label{Sec:Background}

Let $\mathcal{L}_\mathrm{po}$ be the language of partial orders-- first-order logic endowed with a binary (ordering) relation $\sqsubseteq$. We will make reference to the following $\mathcal{L}_\mathrm{po}$ theories:
\begin{itemize}
\item The theory of {\it partial orders} ($\mathrm{PO}$) is the $\mathcal{L}_\mathrm{po}$-theory with axioms asserting that $\sqsubseteq$ is reflexive, antisymmetric and transitive.
\item The theory of {\it lattices} ($\mathrm{Lat}$) is the $\mathcal{L}_\mathrm{po}$-theory extending $\mathrm{PO}$ with axioms asserting that there exists a least element ($0$), and that every pair of elements $x$ and $y$ have both a least upper bound ($x \dot\lor y$) and greatest lower bound ($x \dot\land y$)\footnote{Dots will be used to distinguish the algebraic lattice operations from logical connectives and set-theoretic operations.}.   
\end{itemize}
We use $\mathbf{Atm}(x)$ to abbreviate the $\mathcal{L}_\mathrm{po}$-formula that asserts, in the theory $\mathrm{Lat}$, that $x$ is an atom. I.e. $(x \neq 0) \land \forall y(y \sqsubseteq x \Rightarrow y = 0 \lor y=x)$. 
\begin{itemize}
\item The theory of {\it set-theoretic mereology} ($\mathrm{Mer}$) is the $\mathcal{L}_\mathrm{po}$-theory extending $\mathrm{Lat}$ with the axioms:
\begin{itemize}
\item[](Atomic) for all $x$, there exists $y$ such that $\mathbf{Atm}(y)$ and $y \sqsubseteq x$;
\item[](Unbounded) for all $x$, there exists $y$ such that $x \sqsubseteq y$ and $x \neq y$;
\item[](Relatively Complemented) for all $x$ and $y$, there exists an element $x\dotminus y$ that satisfies the equations
$$y \dot\land (x \dotminus y)= 0 \textrm{ and } x= (x \dot\land y)\dot\lor (x\dotminus y);$$
\item[](Distributive) for all $x$, $y$ and $z$,
$$x \dot\land (y \dot\lor z)= (x \dot\land y) \dot\lor (x \dot\land z)\textrm{ and } x \dot\lor (y \dot\land z)= (x \dot\lor y) \dot\land (x \dot\lor z).$$ 
\end{itemize}
I.e. $\mathrm{Mer}$ is the theory of {\it atomic unbounded relatively complemented distributive lattices}.
\item The theory of {\it infinite atomic Boolean algebras} ($\mathrm{IABA}$) is the $\mathcal{L}_\mathrm{po}$-theory extending $\mathrm{Lat}$ with the axioms:
\begin{itemize}
\item[](Atomic) for all $x$, there exists $y$ such that $\mathbf{Atm}(y)$ and $y \sqsubseteq x$;
\item[](Top) there exists a greatest element ($1$);
\item[](Complemented) for all $x$, there exists an element $\dotminus x$ that satisfies the equations
$$x \dot\lor \dotminus x = 1 \textrm{ and } x \dot\land \dotminus x = 0;$$
\item[](Distributive) for all $x$, $y$ and $z$,
$$x \dot\land (y \dot\lor z)= (x \dot\land y) \dot\lor (x \dot\land z)\textrm{ and } x \dot\lor (y \dot\land z)= (x \dot\lor y) \dot\land (x \dot\lor z);$$
\item[](Infinity Scheme) for all $n \in \mathbb{N}$ with $n > 0$, the axiom that asserts that there are at least $n$ distinct atoms.
\end{itemize}
\end{itemize}

\begin{Definitions1}
Let $\mathcal{M}= \langle M, \sqsubseteq^\mathcal{M} \rangle$ be a model of $\mathrm{IABA}$. A set $I \subseteq M$ is an {\bf proper ideal} of $\mathcal{M}$ if
\begin{itemize}
\item[(i)] $0 \in I$ and $1 \notin I$;
\item[(ii)] for all $x, y \in I$, $x \dot\lor y \in I$;
\item[(iii)] for all $x \in I$ and for all $y \in M$, if $y \sqsubseteq^\mathcal{M} x$, then $y \in I$. 
\end{itemize} 
\end{Definitions1}

Let $\mathcal{L}$ denote the language of set theory-- first-order logic endowed with a binary (membership) relation $\in$. We will have cause to consider the following subsystems of $\mathrm{ZF}$:
\begin{itemize}
\item {\it Adjunctive Set Theory} ($\mathrm{AS}$) is the $\mathcal{L}$-theory with axioms:
\begin{itemize}
\item[]($\mathrm{Emp}$) $\exists x \forall y (y \notin x)$;
\item[]($\mathrm{Adj}$) $\forall x \forall y \exists z \forall u(u \in z \iff (u \in x \lor u=y))$. 
\end{itemize}
The axiom $\mathrm{Adj}$ guarantees that for all sets $x$ and $y$, the set $x \cup \{y\}$ exists.
\item {\it Adjunctive Set Theory with Boolean operations} ($\mathrm{BAS}$) is the $\mathcal{L}$-theory extending $\mathrm{AS}$ with the axioms:
\begin{itemize}
\item[]($\mathrm{Ext}$) $\forall x \forall y(x=y \iff \forall z(z \in x \iff z \in y))$;
\item[]($\mathrm{Union}$) $\forall x \forall y \exists z \forall u(u \in z \iff u \in x \lor u \in y)$;
\item[]($\mathrm{Intersection}$) $\forall x \forall y \exists z \forall u(u \in z \iff u \in x \land u \in y)$;
\item[]($\mathrm{RelComp}$) $\forall x \forall y \exists z \forall u(u \in z \iff u \in x \land u \notin y)$;
\item[]($\mathrm{UB}$) $\forall x \exists y(y \notin x)$.
\end{itemize}
\end{itemize}
Adjunctive Set Theory was first introduced by Tarski, and is known to interpret Robinson's Arithmetic. It follows that both $\mathrm{AS}$ and $\mathrm{BAS}$ are essentially undecidable theories. We refer the reader to \cite{cheXX} for a survey of essentially undecidable theories. As usual, $x \subseteq y$ if $\forall z(z \in x \Rightarrow z \in y)$. If $\phi$ is a formula, then $\mathbf{Var}(\phi)$ will be used to denote the set of variables, both free and bound, that appear in $\phi$.

\begin{Definitions1} 
Let $\mathcal{L}^\prime$ be $\mathcal{L}$ or an extension of $\mathcal{L}$ that is obtained by only adding new unary predicates. Let $\phi$ be an $\mathcal{L}^\prime$-formula and let $n \in \omega$. We say that $\sigma: \mathbf{Var}(\phi) \longrightarrow \mathbb{N}$ is a {\bf stratification} of $\phi$ if 
\begin{itemize}
\item[(i)] if `$x \in y$' is a subformula of $\phi$, then $\sigma(\textrm{`}y\textrm{'})= \sigma(\textrm{`}x\textrm{'})+1$,
\item[(ii)] if `$x=y$' is a subformula of $\phi$, then $\sigma(\textrm{`}y\textrm{'})= \sigma(\textrm{`}x\textrm{'})$. 
\end{itemize}
If there exists a stratification $\sigma$ of $\phi$, then we say that $\phi$ is {\bf stratified}. If there exists a stratification $\sigma: \mathbf{Var}(\phi) \longrightarrow n$ of $\phi$, then we say that $\phi$ is {\bf $n$-stratified}. 
\end{Definitions1}

The notion of stratification was first introduced by Quine in \cite{qui37} where it is used to define the set theory $\mathrm{NF}$. We refer the reader to \cite{for92} for a survey of $\mathrm{NF}$ and related systems of set theory that avoid the set-theoretic paradoxes by restricting comprehension using the notion of stratified formula. Note that the formula $x \subseteq y$ and all of the axioms of $\mathrm{BAS}$ are $2$-stratified $\mathcal{L}$-formulae.

Let $\mathcal{L}_\mathrm{cl}$ be the single-sorted language of class theory-- first-order logic endowed with a binary membership relation ($\in$) and a unary predicate ($\mathcal{S}$) that distinguishes sets from classes. Therefore an $\mathcal{L}_\mathrm{cl}$-structure is a triple $\mathcal{M}= \langle M, \mathcal{S}^\mathcal{M}, \in^\mathcal{M} \rangle$ with $\mathcal{S}^\mathcal{M} \subseteq M$. In order make the presentation of $\mathcal{L}_\mathrm{cl}$-formulae more readable, we will pretend that $\mathcal{L}_\mathrm{cl}$ is a two-sorted language with sorts {\it sets} (referred to using lower case Roman letters $x, y, z, u, v, \ldots$) that are elements of the domain that satisfy $\mathcal{S}$ and {\it classes} (referred to using upper case Roman letters $X, Y, Z, U, V, \ldots$) that are any element of the domain. Therefore, in the axiomatisations presented below, $\exists x(\cdots)$ is an abbreviation for $\exists x (\mathcal{S}(x) \land \cdots)$, $\forall x(\cdots)$ is an abbreviation for $\forall x(\mathcal{S}(x) \Rightarrow \cdots)$ and $\exists x(x=X)$ is an abbreviation for $\mathcal{S}(X)$, etc. We use $\mathrm{NBG}$ to denote the $\mathcal{L}_\mathrm{cl}$-axiomatisation of the version of von Neumann-G\"{o}del-Bernays Class Theory without the axiom of choice that is presented in \cite[Chapter 4]{men97}. We will study the following subsystems of $\mathrm{NBG}$:
\begin{itemize}
\item {\it Adjunctive Set Theory with Classes} ($\mathrm{CAS}$ \footnote{Note that this theory is stronger than the {\it adjunctive class theory} ($\texttt{ac}$) described in \cite[p6-7]{vis10}}) is the $\mathcal{L}_\mathrm{cl}$-theory with axioms:
\begin{itemize}
\item[]($\mathrm{Mem}$) $\forall X \forall Y(X \in Y \Rightarrow \exists x(x= X)$;
\item[]($\mathrm{Subset}$) $\forall x \forall X (\forall y(y \in X \Rightarrow y \in x) \Rightarrow \exists z(z=X))$;
\item[]($\mathrm{Emp}$) $\exists x \forall y (y \notin x)$;
\item[]($\mathrm{Adj}$) $\forall x \forall y \exists z \forall u(u \in z \iff (u \in x \lor u=y))$.
\end{itemize}
\item {\it Adjunctive Class Theory with Boolean operations} ($\mathrm{BAC}$) is the $\mathcal{L}_\mathrm{cl}$-theory extending $\mathrm{CAS}$ with the axioms:
\begin{itemize}
\item[]($\mathrm{CExt}$) $\forall X \forall Y(X=Y \iff \forall x( x \in X \iff x \in Y))$;
\item[]($\mathrm{Union}$) $\forall x \forall y \exists z \forall u(u \in z \iff u \in x \lor u \in y)$;
\item[]($\mathrm{UB}$) $\forall x \exists y(y \notin x)$;
\item[]($\mathrm{CUnion}$) $\forall X \forall Y \exists Z \forall x(x \in Z \iff x \in X \lor x \in Y)$;
\item[]($\mathrm{CIntersection}$) $\forall X \forall Y \exists Z \forall x(x \in Z \iff x \in X \land x \in Y)$;
\item[]($\mathrm{CComp}$) $\forall X \exists Y \forall x(x \in Y \iff x \notin X)$.
\end{itemize} 
\end{itemize} 

Hamkins and Kikuchi \cite[Theorem 9 and Corollary 10]{hk16} observe that the subset relation in any model of $\mathrm{ZF}$ is an atomic unbounded relatively complemented distributive lattice ordering of the universe and, extending \cite{ers64}, show that this theory is complete and decidable.

\begin{Theorems1} \label{Th:HamkinsKikuchi1}
Let $\mathcal{M}= \langle M, \in^\mathcal{M} \rangle$ be a model of $\mathrm{ZF}$. Then\\ \mbox{$\langle M, \subseteq^\mathcal{M} \rangle \models \mathrm{Mer}$} and this theory is complete and decidable. \Square 
\end{Theorems1}

\cite{hk16} note that this result also holds for models of certain subsystems of $\mathrm{ZF}$ such as finite set theory and Kripke-Platek Set Theory. In the next section we will see that Theorem \ref{Th:HamkinsKikuchi1} holds when $\mathcal{M}$ satisfies $\mathrm{BAS}$. 

In \cite[Section 5]{hkxx}, Hamkins and Kikuchi extend their analysis of the theory and structure of the subset ordering to models of class theory. They show that if $\mathcal{M}$ is a model of Von Neumann-Bernays-G\"{o}del set theory with a class version of the axiom of choice, then the structure $\langle M, \subseteq^\mathcal{M} \rangle$ is an $\omega$-saturated model of $\mathrm{IABA}$. Even if $\mathcal{M}$ does not satisfy the axiom of choice, $\langle M, \subseteq^\mathcal{M} \rangle$ satisfies $\mathrm{IABA}$. It is a well-known result due to Tarski \cite{tar49} that the theory $\mathrm{IABA}$ complete and decidable.  

\begin{Theorems1} \label{Th:IABAComplete}
Let $\mathcal{M}= \langle M, \mathcal{S}^\mathcal{M}, \in^\mathcal{M} \rangle$ be a model of $\mathrm{NBG}$. Then $\langle M, \subseteq^\mathcal{M} \rangle \models \mathrm{IABA}$ and this theory is complete and decidable. \Square  
\end{Theorems1}

\section[The set theory corresponding to $\mathrm{Mer}$]{The set theory corresponding to $\mathrm{Mer}$} \label{Sec:SetTheoryAndMer}

This section establishes that $\mathrm{ZF}$ can be replaced by $\mathrm{BAS}$ in Theorem \ref{Th:HamkinsKikuchi1} and, conversely, that every model of $\mathrm{Mer}$ with the same number of atoms as elements can be realised as the subset relation of a model of $\mathrm{BAS}$.

\begin{Theorems1} \label{Th:MainTheorem1}
Let $\mathcal{M}=\langle M, \in^M \rangle$ be a model of $\mathrm{BAS}$. Then $\langle M, \subseteq^\mathcal{M} \rangle$ is a model of $\mathrm{Mer}$ with the same number of atoms as elements. 
\end{Theorems1}

\begin{proof}
Work inside $\mathcal{M}$. Now, $\langle M, \subseteq^{\mathcal{M}}\rangle$ is clearly reflexive and transitive, and $\mathrm{Ext}$ ensures that $\langle M, \subseteq^{\mathcal{M}}\rangle$ is antisymmetric. Now, the axiom $\mathrm{Emp}$ ensures that $\emptyset$ exists and this set is a $\subseteq^{\mathcal{M}}$-least element. The axioms $\mathrm{Union}$ and $\mathrm{Intersection}$ ensure that for all sets $x$ and $y$, $x \cap y$ and $x \cup y$ exist. In the order $\langle M, \subseteq^\mathcal{M} \rangle$, $x \cap y$ is the greatest lower bound of $x$ and $y$, and $x \cup y$ is the least upper bound of $x$ and $y$. Since for all $x$, $y$ and $z$,
$$x \cap (y \cup z)= (x \cap y)\cup (x\cap z) \textrm{ and } x \cup (y \cap z)= (x \cup y)\cap (x\cup z),$$
$\langle M, \subseteq^{\mathcal{M}}\rangle$ is a distributive lattice.
The existence of $\emptyset$ and the axiom $\mathrm{Adj}$ ensures that for all $x$, $\{x\}$ exists. Note that for all $x$, $\{x\}$ is an atom of $\langle M, \subseteq^\mathcal{M} \rangle$. Now, if $y \neq \emptyset$, then there exists $x \in y$ and $\{x\}$ is an atom below $y$ in $\langle M, \subseteq^\mathcal{M} \rangle$. Therefore, $\langle M, \subseteq^\mathcal{M} \rangle$ is atomic. The axiom $\mathrm{RelComp}$ ensures that for all sets $x$ and $y$, $x-y$ exists. Since for all $x$ and $y$,
$$y \cap (x-y)= \emptyset \textrm{ and } x= (x \cap y) \cup (x-y),$$
$\langle M, \subseteq^\mathcal{M} \rangle$ is relatively complemented. Finally, for all $x$, the axiom $\mathrm{UB}$ ensures that there exists $y$ such that $y \notin x$. The axiom $\mathrm{Adj}$ ensures that $x \cup \{y\}$ exists. Now, $x \subseteq x \cup \{y\} \neq x$ and so $\langle M, \subseteq^\mathcal{M} \rangle$ is unbounded. To see that $\langle M, \subseteq^\mathcal{M} \rangle$ has the same number of atoms as elements, observe that $f= \{\langle x, \{x\} \rangle \mid x \in M \}$ is a bijection between the elements of $\langle M, \subseteq^\mathcal{M} \rangle$ and the atoms of $\langle M, \subseteq^\mathcal{M} \rangle$.   
\Square    
\end{proof} 
   
We now turn to showing that if $\mathcal{M}= \langle M, \sqsubseteq^{\mathcal{M}} \rangle$ is a model of $\mathrm{Mer}$ with the same number of atoms as elements, then there exists a membership relation on $M$ that makes $M$ a model of $\mathrm{BAS}$ whose subset relation is exactly $\sqsubseteq^{\mathcal{M}}$. The intuition behind this result is that a model of set theory can be obtained from a subset relation and a map that sends a set to its own singleton (see \cite[Theorem 13]{hk16}). The set theory $\mathrm{BAS}$ is so weak that this singleton map can be chosen to be any bijection between elements and atoms of a structure satisfying $\mathrm{Mer}$. 

\begin{Lemma1} \label{Th:OrderDeterminedByAtoms}
Let $\mathcal{M}= \langle M, \sqsubseteq^{\mathcal{M}} \rangle$ be a model of $\mathrm{Mer}$. For all $x, y \in M$, $x \sqsubseteq^{\mathcal{M}} y$ if and only if for all atoms $z \in M$, if $z \sqsubseteq^{\mathcal{M}} x$, then $z \sqsubseteq^{\mathcal{M}} y$.
\end{Lemma1}

\begin{proof}
Let $x, y \in M$. Work inside $\mathcal{M}$. Note that if $x \sqsubseteq y$, then for all atoms $z$, if $z \sqsubseteq x$, then $z \sqsubseteq y$. Conversely, suppose that for all atoms $z$, if $z \sqsubseteq x$, then $z \sqsubseteq y$. Now, $x \dotminus y =0$, because otherwise there would be an atom that sits below both $x \dotminus y$ and $y$ contradicting the fact that $y \dot\land (x \dotminus y)= 0$. Therefore, $x= x \dot\land y$ and so $x \sqsubseteq y$.   
\Square
\end{proof}

\begin{Theorems1} \label{Th:MainTheorem2}
Let $\mathcal{M}= \langle M, \sqsubseteq^\mathcal{M} \rangle$ be a model of $\mathrm{Mer}$. Let \mbox{$A= \{x \in M \mid \mathcal{M} \models \mathbf{Atm}(x)\}$} and let $f: M \longrightarrow A$ be a bijection. Define $\in^* \subseteq M \times M$ by: for all $x, y \in M$,
$$x \in^* y \textrm{ if and only if } \mathcal{M} \models f(x) \sqsubseteq y.$$
Then $\langle M, \in^* \rangle \models \mathrm{BAS}$ and $\subseteq^{\langle M, \in^*\rangle}= \sqsubseteq^\mathcal{M}$.  
\end{Theorems1}

\begin{proof}
It follows immediately from Lemma \ref{Th:OrderDeterminedByAtoms} that $\langle M, \in^* \rangle \models \mathrm{Ext}$ and $\subseteq^{\langle M, \in^*\rangle}= \sqsubseteq^\mathcal{M}$. Note that the $\sqsubseteq^\mathcal{M}$-least element $0$ is such that $\langle M, \in^* \rangle \models \forall y(y \notin 0)$. Therefore \mbox{$\langle M, \in^* \rangle \models \mathrm{Emp}$}. Now, let $x, y \in M$. Working in $\mathcal{M}$, let $z= x \dot\lor f(y)$. For all $u \in M$,
$$f(u) \sqsubseteq^\mathcal{M} z \textrm{ if and only if } f(u) \sqsubseteq^\mathcal{M} x \textrm{ or } f(u)=f(y)$$ 
$$\textrm{if and only if } \langle M, \in^* \rangle \models (u \in x \lor u=y).$$
Therefore $\langle M, \in^* \rangle \models \mathrm{Adj}$. To show that Union holds in $\langle M, \in^* \rangle$, let $x, y \in M$. Working in $\mathcal{M}$, let $z= x \dot\lor y$. Let $u \in M$. It is clear that if $f(u) \sqsubseteq^{\mathcal{M}} x$ or $f(u) \sqsubseteq^{\mathcal{M}} y$, then $f(u) \sqsubseteq^{\mathcal{M}} z$. Conversely, suppose that $f(u) \not\sqsubseteq^{\mathcal{M}} x$ and $f(u) \not\sqsubseteq^{\mathcal{M}} y$. Then, working inside $\mathcal{M}$, $x \dot\land f(u)= y \dot\land f(u)= 0$ and 
$$f(u) \dot\land z= f(u)\dot\land (x \dot\lor y)= (f(u) \dot\land x)\dot\lor (f(u) \dot\land y)= 0.$$
Therefore $f(u) \sqsubseteq^{\mathcal{M}} z$. We have shown that 
$$\langle M, \in^* \rangle \models \forall u( u \in z \iff u \in x \lor u \in y).$$
So, $\langle M, \in^* \rangle \models \mathrm{Union}$. To see that $\langle M, \in^*\rangle$ satisfies $\mathrm{Intersection}$, let $x, y \in M$. Working inside $\mathcal{M}$, let $z= x \dot\land y$. Let $u \in M$. If $f(u) \sqsubseteq^{\mathcal{M}} z$, then $f(u) \sqsubseteq^{\mathcal{M}} x$ and $f(u) \sqsubseteq^{\mathcal{M}} y$. Conversely, suppose $f(u) \sqsubseteq^{\mathcal{M}} x$ and $f(u) \sqsubseteq^{\mathcal{M}} y$. Therefore, inside $\mathcal{M}$, 
$$f(u) \dot\lor z= f(u) \dot\lor (x \dot\land y)= (f(u) \dot\lor x) \dot\land (f(u) \dot\lor y)= x \dot\land y= z.$$
Therefore $f(u) \sqsubseteq^{\mathcal{M}} z$. This shows that
$$\langle M, \in^* \rangle \models \forall u (u \in z \iff u \in x \land u \in y).$$
So, $\langle M, \in^*\rangle \models \mathrm{Intersection}$. To show that $\langle M, \in^* \rangle$ satisfies $\mathrm{RelComp}$, let $x, y \in M$. Inside $\mathcal{M}$, let $z= x \dotminus y$. Let $u \in M$. Note that, inside $\mathcal{M}$, $x \dotminus y \sqsubseteq x$ and, since $y \dot\land (x \dotminus y)= 0$, if $f(u) \sqsubseteq x \dotminus y$, then $f(u) \not\sqsubseteq y$. Therefore, if $f(u) \sqsubseteq^{\mathcal{M}} z$, then $f(u) \sqsubseteq^{\mathcal{M}} x$ and $f(u) \not\sqsubseteq^{\mathcal{M}} y$. Conversely, suppose that $f(u) \sqsubseteq^{\mathcal{M}} x$ and $f(u) \not\sqsubseteq^{\mathcal{M}} y$. Therefore, inside $\mathcal{M}$,
$$f(u)= f(u) \dot\land x = f(u) \dot\land ((x \dot\land y) \dot\lor z)= (f(u) \dot\land (x \dot\land y)) \dot\lor (f(u) \dot\land z)= 0 \dot\lor (f(u) \dot \land z)= f(u) \dot\land z.$$
Therefore, $f(u) \sqsubseteq^{\mathcal{M}} z$. This shows that 
$$\langle M, \in^* \rangle \models \forall u( u \in z \iff u \in x \land u \notin y).$$
So, $\langle M, \in^* \rangle \models \mathrm{RelComp}$. To see that $\langle M, \in^* \rangle$ satisfies $\mathrm{UB}$, let $x \in M$. Since $\mathcal{M}$ is unbounded, there exists $z \in M$ such that $x \sqsubseteq^\mathcal{M} z$ and $x \neq z$. Since $x \neq z$, by Lemma \ref{Th:OrderDeterminedByAtoms}, there is an atom $u$ such $u \sqsubseteq^{\mathcal{M}} z$ and $u \not\sqsubseteq^{\mathcal{M}} x$. Now, letting $y= f^{-1}(u)$, we see that $\langle M, \in^*\rangle \models (y \notin x)$. Therefore $\langle M, \in^*\rangle \models \mathrm{UB}$.
\Square         
\end{proof}  

\begin{Coroll1} \label{Th:CountableModelsOfMerAreModelsOfBAS}
Let $\mathcal{M}= \langle M, \sqsubseteq^{\mathcal{M}} \rangle$ be a countable model of $\mathrm{Mer}$. There exists $\in^* \subseteq M \times M$ such that $\langle M, \in^* \rangle \models \mathrm{BAS}$ and $\subseteq^{\langle M, \in^*\rangle}= \sqsubseteq^{\mathcal{M}}$.
\end{Coroll1}

\begin{proof}
Let $A= \{x \in M \mid \mathcal{M} \models \mathbf{Atm}(x)\}$. Note that $A$ is infinite. Since $\mathcal{M}$ is countable, $|M|=|A|$ and the result now follows from Theorem \ref{Th:MainTheorem2}.
\Square
\end{proof}

In contrast with Corollary \ref{Th:CountableModelsOfMerAreModelsOfBAS}, there are uncountable models of $\mathrm{Mer}$ with fewer atoms than elements and thus models of $\mathrm{Mer}$ that are not isomorphic to the subset relation of any model of $\mathrm{BAS}$. 

\begin{Examp1} \label{Ex:ModelOfMerThatCannotBeSubset}
Let $X= \{ B \subseteq \omega+\omega \mid (\exists \beta \in \omega+\omega)(\forall \alpha \in B)(\alpha \in \beta)\}$. Note that $\langle X, \subseteq \rangle$ is a model of $\mathrm{Mer}$. The atoms of $\langle X, \subseteq \rangle$ are
$$A= \{\{\alpha\} \mid \alpha \in \omega+\omega\}.$$
Since $|A|= \aleph_0 < |X|= 2^{\aleph_0}$, Theorem \ref{Th:MainTheorem1} shows that $\langle X, \subseteq \rangle$ is not isomorphic to the subset relation of any model of $\mathrm{BAS}$.  
\end{Examp1}

Theroem \ref{Th:MainTheorem2} also shows that $\mathrm{BAS}$ is the minimal set theory whose subset relation is an atomic unbounded relatively complemented lattice order.

\begin{Theorems1} \label{Th:BASMinimal}
Let $\mathcal{M}= \langle M, \in^\mathcal{M} \rangle$ is a model of $\mathrm{AS}$. Then the following are equivalent:
\begin{itemize}
\item[(I)] $\mathcal{M} \models \mathrm{BAS}$;
\item[(II)] $\langle M, \subseteq^\mathcal{M} \rangle \models \mathrm{Mer}$.
\end{itemize}
\end{Theorems1}

\begin{proof}
The fact that $(I) \Rightarrow (II)$ follows from \ref{Th:MainTheorem1}. To see that $(II) \Rightarrow (I)$, assume that $\langle M, \subseteq^\mathcal{M} \rangle \models \mathrm{Mer}$. The fact that $\langle M, \subseteq^\mathcal{M} \rangle$ is antisymmetric immediately implies that $\mathcal{M} \models \mathrm{Ext}$. Inside $\mathcal{M}$, the axioms $\mathrm{Emp}$ and $\mathrm{Adj}$ ensure that for all sets $x$, $\{x\}$ exists. Consider $f= \{ \langle x, \{x\} \rangle \mid x \in M\}$. The existence of singletons coupled with extensionality in $\mathcal{M}$ ensure that $f$ is a bijection between $M$ and the atoms of $\langle M, \subseteq^\mathcal{M} \rangle$. As in Theorem \ref{Th:MainTheorem2}, define $\in^* \subseteq M \times M$ by: for all $x, y \in M$,
$$x \in^* y \textrm{ if and only if } \mathcal{M} \models f(x) \subseteq^\mathcal{M} y.$$
Therefore $\langle M, \in^*\rangle \models \mathrm{BAS}$. But, for all $x, y \in M$,
$$x \in^* y \textrm{ if and only if } \{x\} \subseteq^{\mathcal{M}} y \textrm{ if and only if } x \in^\mathcal{M} y.$$
Therefore $\langle M, \in^\mathcal{M} \rangle \models \mathrm{BAS}$. 
\Square
\end{proof}

\section[The subset relation and $2$-stratified sentences]{The subset relation and $2$-stratified sentences} \label{Sec:2StratifiedSentences}

In this section we show every $2$-stratified $\mathcal{L}$-sentence, $\phi$, can be translated into an $\mathcal{L}_\mathrm{po}$-sentence, $\phi^\tau$, such that if $\mathcal{M}= \langle M, \in^\mathcal{M} \rangle$ is a model of $\mathrm{BAS}$, then $\phi$ holds in $\mathcal{M}$ if and only if $\phi^\tau$ holds in $\langle M, \subseteq^\mathcal{M} \rangle$. Since the axioms of $\mathrm{BAS}$ completely determine the theory of $\langle M, \subseteq^\mathcal{M} \rangle$, this shows that $\mathrm{BAS}$ decides every $2$-stratified sentence.

We begin with the observation that by replacing the order $\sqsubseteq$ by the set-theoretic $\subseteq$ order in an $\mathcal{L}_\mathrm{po}$-formula one obtains a $2$-stratified $\mathcal{L}$-formula.

\begin{Theorems1}
Let $\phi(x_1, \ldots, x_n)$ be an $\mathcal{L}_\mathrm{po}$-formula. The $\mathcal{L}$-formula $\phi^S(\vec{x})$ obtained by replacing $\sqsubseteq$ by the set-theoretic definition of $\subseteq$ admits a stratification $\sigma: \mathbf{Var}(\phi^S) \longrightarrow 2$ such that for all $1 \leq i \leq n$, $\sigma(\textrm{`}x_i\textrm{'})= 1$.
\end{Theorems1}

\begin{proof}
The $\mathcal{L}$-formulae $x=y$ and $x \subseteq y$ admit $2$-stratifications that assigned the value $1$ to both the variables $x$ and $y$. The theorem now follows by a straightforward induction on the structural complexity of $\phi(x_1, \ldots, x_n)$. \Square
\end{proof}

We now turn to defining a translation of $2$-stratified $\mathcal{L}$-formulae into $\mathcal{L}_\mathrm{po}$-formulae.

\begin{Definitions1} \label{Df:TranslationOf2StratFormulae}
Let $\phi$ be an $\mathcal{L}$-formula with a stratification $\sigma: \mathbf{Var}(\phi) \longrightarrow 2$. Define a translation, $\tau$, of $\phi$ into an $\mathcal{L}_\mathrm{po}$-formula $\phi^\tau$ by recursion:
$$\textrm{if } \phi \textrm{ is } u \in v, \textrm{ then } \phi^\tau \textrm{ is } \mathbf{Atm}(u) \land u \sqsubseteq v;$$
$$\textrm{if } \phi \textrm{ is } u=v, \textrm{ then } \phi^\tau \textrm{ is } u=v;$$
$$\textrm{if } \phi \textrm{ is } \psi \# \theta, \textrm{ where } \# \in \{\lor, \land, \Rightarrow, \iff\}, \textrm{ then } \phi^\tau \textrm{ is } \psi^\tau \# \theta^\tau;$$
$$\textrm{if } \phi \textrm{ is } \neg \psi, \textrm{ then } \phi^\tau \textrm{ is } \neg \psi^\tau;$$
$$\textrm{if } \phi \textrm{ is }\exists v\ \psi, \textrm{ then } \phi^\tau \textrm{ is } \left\{\begin{array}{ll}
\exists v (\mathbf{Atm}(v) \land \psi) & \textrm{ if } \sigma(\textrm{`}v\textrm{'})= 0\\
\exists v\ \psi & \textrm{ if }\sigma(\textrm{`}v\textrm{'})= 1 
\end{array}\right.;$$
$$\textrm{if } \phi \textrm{ is }\forall v\ \psi, \textrm{ then } \phi^\tau \textrm{ is } \left\{\begin{array}{ll}
\forall v (\mathbf{Atm}(v) \Rightarrow \psi) & \textrm{ if } \sigma(\textrm{`}v\textrm{'})= 0\\
\forall v\ \psi & \textrm{ if }\sigma(\textrm{`}v\textrm{'})= 1 
\end{array}\right.$$
\end{Definitions1}

Note that $\tau$ treats variables that appear on the left-hand side of the membership relation differently from variables that appear on the right-hand side of the membership relation. This prevents this translation from being generalised to formulae which contain two atomic subformulae that contain the same variable appearing on different sides of the membership relation.   

\begin{Theorems1} \label{Th:TranslationOf2StratToSubsetAlgebra}
Let $\phi(x_1, \ldots, x_n, y_1, \ldots, y_m)$ be an $\mathcal{L}$-formula with a stratification $\sigma: \mathbf{Var}(\phi) \longrightarrow 2$ such that for all $1 \leq i \leq n$, $\sigma(\textrm{`}x_i\textrm{'})= 0$, and for all $1 \leq i \leq m$, $\sigma(\textrm{`}y_i\textrm{'})= 1$. Let $\mathcal{M}= \langle M, \in^\mathcal{M} \rangle$ be a model of $\mathrm{BAS}$. For all $a_1, \ldots, a_n, b_1, \ldots, b_m \in M$,
$$\mathcal{M} \models \phi(a_1, \ldots, a_n, b_1, \ldots, b_m) \textrm{ if and only if}$$
$$\langle M, \subseteq^\mathcal{M} \rangle \models \phi^\tau(\{a_1\}, \ldots, \{a_n\}, b_1, \ldots, b_m).$$  
\end{Theorems1}

\begin{proof}
For all $a, b \in M$, 
$$\mathcal{M} \models a \in b \textrm{ if and only if } \langle M, \subseteq^\mathcal{M} \rangle \models \{a\} \subseteq b.$$
The theorem then follows by a straightforward induction on the structural complexity of $\phi(x_1, \ldots, x_n, y_1, \ldots, y_m)$.
\Square
\end{proof}

Combined with the fact that the theory of the subset relation the set theory $\mathrm{BAS}$ is a complete theory, this result shows that $\mathrm{BAS}$ decides every $2$-stratified sentence. 

\begin{Coroll1}
If $\phi$ is a $2$-stratified $\mathcal{L}$-sentence, then 
$$\mathrm{BAS} \vdash \phi \textrm{ or } \mathrm{BAS} \vdash \neg \phi.$$
\end{Coroll1}

\begin{proof}
Let $\phi$ be a $2$-stratified $\mathcal{L}$-sentence. Since $\mathrm{Mer}$ is complete \ref{Th:HamkinsKikuchi1},
$$\mathrm{Mer} \vdash \phi^\tau \textrm{ or } \mathrm{Mer} \vdash \neg \phi^\tau.$$
It follows from Theorems \ref{Th:MainTheorem1} and \ref{Th:TranslationOf2StratToSubsetAlgebra} that either $\phi$ holds in all models of $\mathrm{BAS}$ or $\neg \phi$ holds in all models of $\mathrm{BAS}$. 
\Square
\end{proof}

Gogol \cite{gog78} shows that $\mathrm{ZFC}$ decides every sentence in prenex normal form that is prefixed by a block of universal quantifiers followed by a single existential quantifier. H. Friedman \cite{fri03} extending the work of Gogol \cite{gog79} shows that a weak subsystem of $\mathrm{ZF}$ decides every sentence with only three quantifiers. 

\begin{Coroll1}
Let $T \supseteq \mathrm{BAS}$ be consistent. Then 
$$\{\ulcorner \phi \urcorner \mid (\phi \textrm{ is a }2\textrm{-stratified }\mathcal{L}\textrm{-sentence})\land (T \vdash \phi)\} \textrm{ is recursive}.$$
\end{Coroll1}

\begin{proof}
The algorithm that tests whether $\phi$ is a $2$-stratified $\mathcal{L}$-sentence and then simultaneously searches for a proof of $\phi$ or $\neg \phi$ from the axioms of $\mathrm{BAS}$ decides this set and will always halt.
\Square 
\end{proof}

In particular, $\mathrm{ZF}$ decides every $2$-stratified sentence. We conclude this section by observing that this result cannot be extended to $3$-stratified sentences.

\begin{Theorems1}
There is a $3$-stratified $\mathcal{L}$-sentence that is not decided by $\mathrm{ZF}$.
\end{Theorems1}

\begin{proof}
The following $\mathcal{L}$-sentence asserts that {\it for every non-empty set of non-empty disjoint sets $X$, there exists a set $C$ such that for all $x \in X$, $x \cap C$ is a singleton}:
$$\forall X \left(\begin{array}{c}
(\forall x \in X) \exists w (w \in x) \land (\forall x, y \in X)(x \neq y \Rightarrow \neg \exists w (w \in x \land w \in y)) \Rightarrow\\
\exists C \left(\begin{array}{c}
(\forall z \in C)(\exists x \in X)(z \in x) \land (\forall x \in X)(\exists z \in C)(z \in x)\\
\land (\forall z, w \in C)(\forall x \in X)(z \in x \land w \in x \Rightarrow z=w)
\end{array}  
\right)
\end{array}\right).$$
This sentence is $3$-stratified and equivalent to the Axiom of Choice. 
\Square
\end{proof}

\section[The class theory corresponding to $\mathrm{IABA}$]{The class theory corresponding to $\mathrm{IABA}$} \label{Sec:ClassTheoryAndIABA}

This section extends the results of sections \ref{Sec:SetTheoryAndMer} and \ref{Sec:2StratifiedSentences} to class theory. We begin by showing that $\mathrm{BAC}$ is the minimal subsystem of $\mathrm{NBG}$ guaranteeing that the definable subset relation is an infinite atomic Boolean algebra.

\begin{Theorems1} \label{Th:SubsetRelationInBAC}
Let $\mathcal{M}= \langle M, \mathcal{S}^\mathcal{M}, \in^\mathcal{M} \rangle$ be a model of $\mathrm{BAC}$. Then $\langle M, \subseteq^\mathcal{M} \rangle$ is a model of $\mathrm{IABA}$, $\mathcal{S}^\mathcal{M}$ is a proper ideal of $\langle M, \subseteq^\mathcal{M} \rangle$ that contains all of the atoms of $\langle M, \subseteq^\mathcal{M} \rangle$ and $\mathcal{S}^\mathcal{M}$ is the same size as the set of atoms of $\langle M, \subseteq^\mathcal{M} \rangle$.    
\end{Theorems1}

\begin{proof}
Work inside $\mathcal{M}$. Since $\mathcal{M}$ satisfies $\mathrm{Subset}$, $\mathcal{S}^\mathcal{M}$ is closed downwards by $\subseteq^\mathcal{M}$ in $\langle M, \subseteq^\mathcal{M} \rangle$. Since $\langle \mathcal{S}^\mathcal{M}, \in^\mathcal{M} \rangle \models \mathrm{BAS}$, it follows immediately from Theorem \ref{Th:MainTheorem1} that $\langle \mathcal{S}^\mathcal{M}, \subseteq^\mathcal{M} \rangle$ is a model of $\mathrm{Mer}$ with the same number of atoms as elements. Therefore, $\mathcal{S}^\mathcal{M}$ is a proper ideal of $\langle M, \subseteq^\mathcal{M} \rangle$ with the same number of atoms as elements. It is clear that the atoms of $\langle M, \subseteq^\mathcal{M} \rangle$ all belong to $\mathcal{S}^\mathcal{M}$. We are left to show that $\langle M, \subseteq^\mathcal{M} \rangle$ is an infinite atomic Boolean algebra. Note that $\langle M, \subseteq^\mathcal{M} \rangle$ is reflexive and transitive, and $\mathrm{CExt}$ ensures that $\langle M, \subseteq^\mathcal{M} \rangle$ is antisymmetric. The axiom $\mathrm{Emp}$ ensures that $\emptyset$ exists and this set is the least element of $\langle M, \subseteq^\mathcal{M} \rangle$. The axioms $\mathrm{Emp}$ and $\mathrm{CComp}$ ensure that there exists a class, $V$, that contains every set, and this class is the greatest element of $\langle M, \subseteq^\mathcal{M} \rangle$. In the order $\langle M, \subseteq^\mathcal{M} \rangle$, $X \cup Y$ is the least upper bound of $X$ and $Y$, and $X \cap Y$ is the greatest lower bound of $X$ and $Y$. Since for all $X$, $Y$ and $Z$,
$$X \cap (Y \cup Z)= (X \cap Y) \cup (X \cap Z) \textrm{ and } X \cup (Y \cap Z)= (X \cup Y) \cap (X \cup Z),$$
$\langle M, \subseteq^\mathcal{M} \rangle$ is a distributive lattice. If $X$ is a class and $x$ is a set with $x \in X$, then $\{x\}$ is an atom that sits below $X$ in $\langle M, \subseteq^\mathcal{M} \rangle$. Therefore, $\langle M, \subseteq^\mathcal{M} \rangle$ is atomic. Let $X$ be a class. The axiom $\mathrm{CComp}$ ensures that $V - X$ is a class and 
$$V= X \cup (V-X) \textrm{ and } \emptyset= X \cap (V-X).$$
Therefore, $\langle M, \subseteq^\mathcal{M} \rangle$ has complements. The fact that $\mathcal{S}^\mathcal{M}$ is a proper ideal of $\langle M, \subseteq^\mathcal{M} \rangle$ shows that $\langle M, \subseteq^\mathcal{M} \rangle$ has infinitely many atoms.
\Square
\end{proof}

We now turn to showing that every infinite atomic Boolean algebra can be realised as the subset relation of a model of $\mathrm{BAC}$. Lemma \ref{Th:OrderDeterminedByAtoms} also holds for infinite atomic Boolean algebras.

\begin{Lemma1} \label{Th:ElementsDeterminedByAtomsIABA}
Let $\mathcal{M}= \langle M, \sqsubseteq^{\mathcal{M}} \rangle$ be a model of $\mathrm{IABA}$. For all $x, y \in M$, $x \sqsubseteq^{\mathcal{M}} y$ if and only if for all atoms $z \in M$, if $z \sqsubseteq^{\mathcal{M}} x$, then $z \sqsubseteq^{\mathcal{M}} y$.
\end{Lemma1}

\begin{proof}
Identical to the proof of Lemma \ref{Th:OrderDeterminedByAtoms}.\Square
\end{proof}

\begin{Lemma1} \label{Th:ProperIdealsContainingAtomsSatisfyMer}
Let $\mathcal{M}= \langle M, \sqsubseteq^{\mathcal{M}} \rangle$ be a model of $\mathrm{IABA}$. If $I \subseteq M$ is a proper ideal of $\mathcal{M}$ that contains all of the atoms of $\mathcal{M}$, then $\mathcal{I}= \langle I, \sqsubseteq^\mathcal{M} \rangle$ is a model of $\mathrm{Mer}$ such that for all $x, y \in I$,
$$(x \dot\lor y)^\mathcal{I}= (x \dot\lor y)^\mathcal{M},\ (x \dot\land y)^\mathcal{I}= (x \dot\land y)^\mathcal{M} \textrm{ and } (x\dotminus y)^\mathcal{I}= (x \land \dotminus y)^\mathcal{M}.$$  
\end{Lemma1}

\begin{proof}
Let $I \subseteq M$ be a proper ideal of $\mathcal{M}$ that contains all of the atoms of $\mathcal{M}$ and let $\mathcal{I}= \langle I, \sqsubseteq^\mathcal{M} \rangle$. Since $0 \in I$ and $\mathcal{I}$ is closed under least upper bounds and greatest lower bounds from $\mathcal{M}$, $\mathcal{I}$ is a sub-lattice of $\mathcal{M}$ that is distributive, and the operations $\dot\lor$ and $\dot\land$ are the same in $\mathcal{I}$ as they are in $\mathcal{M}$. Since $\mathcal{I}$ contains all of the atoms from $\mathcal{M}$, $\mathcal{I}$ is atomic. To see that $\mathcal{I}$ has relative complements, let $x, y \in I$. Now, $x \dot\land \dotminus y \sqsubseteq^\mathcal{M} x$, so $x \dot\land \dotminus y \in I$. And, in $\mathcal{M}$,  
$$y \dot\land (x \dot\land \dotminus y)= x \dot\land (y \dot\land \dotminus y)= 0 \textrm{ and}$$
$$(x \dot\land y) \dot\lor (x \dot\land \dotminus y)= ((x \dot\land y) \dot \lor x) \dot\land ((x \dot\land y) \dot\lor \dotminus y)= x \dot\land (x \dot\lor \dotminus y)= x,$$
so $\mathcal{I}$ is relatively complemented. Finally, to see that $\mathcal{I}$ is unbounded, let $x \in I$. Since $I$ is a proper ideal and by Lemma \ref{Th:ElementsDeterminedByAtomsIABA}, there exists an atom $a \in I$ such that $a \not \sqsubseteq^\mathcal{M} x$. Now, $x \sqsubseteq^\mathcal{M} x \dot\lor a$ and $x \neq x \dot\lor a$, so $\mathcal{I}$ is unbounded.       
\Square
\end{proof}

These results allow us to show that a model of $\mathrm{BAC}$ can be built from an infinite atomic Boolean algebra equipped with a proper ideal containing all of the atoms that is the same cardinality as the set of atoms.

\begin{Theorems1} \label{Th:ConstructionOfModelsOfBAC}
Let $\mathcal{M}= \langle M, \sqsubseteq^\mathcal{M} \rangle$ be a model of $\mathrm{IABA}$. Let $A= \{x \in M \mid \mathcal{M} \models \mathbf{Atm}(x) \}$. Let $I \subseteq M$ be a proper ideal of $\mathcal{M}$ with $A \subseteq I$ and let $f: I \longrightarrow A$ be a bijection. Define $\in^* \subseteq M \times M$ by: for all $x, y \in M$,
$$x \in^* y \textrm{ if and only if } \mathcal{M} \models f(x) \sqsubseteq y.$$
Then $\langle M, I, \in^* \rangle \models \mathrm{BAC}$ and $\subseteq^{\langle M, I, \in^*\rangle}= \sqsubseteq^\mathcal{M}$. 
\end{Theorems1}

\begin{proof}
Lemma \ref{Th:ElementsDeterminedByAtomsIABA} immediately implies that $\langle M, I, \in^* \rangle \models \mathrm{CExt}$ and $\subseteq^{\langle M, I, \in^*\rangle}= \sqsubseteq^\mathcal{M}$. Lemma \ref{Th:ProperIdealsContainingAtomsSatisfyMer} and Theorem \ref{Th:MainTheorem2} show that $\langle I, \in^* \rangle \models \mathrm{BAS}$. Therefore, $\mathrm{Union}$, $\mathrm{Adj}$, $\mathrm{Emp}$ and $\mathrm{UB}$ hold in $\langle M, I, \in^* \rangle$. The fact that the domain of $f$ is $I$ ensures that $\mathrm{Mem}$ holds in $\langle M, I, \in^* \rangle$. Since $I$ is downward closed, $\langle M, I, \in^* \rangle$ satisfies $\mathrm{Subset}$. If $x, y \in M$, then the same arguments used in the proof of Theorem \ref{Th:MainTheorem2} show that the objects $x \dot\land y$ and $x \dot\lor y$ in $\mathcal{M}$ witness the fact that $\mathrm{CIntersection}$ and $\mathrm{CUnion}$, respectively, hold for $x$ and $y$. Finally, to verify that $\mathrm{CComp}$ holds, let $x \in M$. Since $x \dot\land \dotminus x= 0$ and $x \dot\lor \dotminus x= 1$ in $\mathcal{M}$, for all atoms $u \in M$, $u \sqsubseteq^\mathcal{M} x$ if and only if $u \not\sqsubseteq \dotminus x$. This shows that $(\dotminus x)^\mathcal{M}$ is the complement of $x$ in $\langle M, I, \in^* \rangle$.    
\Square
\end{proof}   

In contrast to Example \ref{Ex:ModelOfMerThatCannotBeSubset}, every infinite atomic Boolean algebra can be realised as the subset relation of a model of $\mathrm{BAC}$. 

\begin{Theorems1} \label{Th:EveryIABAIsASubsetRelation}
Let $\mathcal{M}= \langle M, \sqsubseteq^\mathcal{M} \rangle$ be a model of $\mathrm{IABA}$. Then there exists $\mathcal{N}= \langle N, \mathcal{S}^\mathcal{N}, \in^\mathcal{N}\rangle$ such that $\mathcal{N} \models \mathrm{BAC}$ and $\langle N, \subseteq^\mathcal{N} \rangle \cong \mathcal{M}$.  
\end{Theorems1}

\begin{proof}
Let $A= \{x \in M \mid \mathcal{M} \models \mathbf{Atm}(x)\}$. Let $I \subseteq M$ be the proper ideal generated by the atoms of $\mathcal{M}$. I.e. 
$$I= \{x \in M \mid \textrm{there are finitely many atoms below } x \textrm{ in } \mathcal{M}\}.$$
Now, by Lemma \ref{Th:ElementsDeterminedByAtomsIABA}, 
$$|I|= |[A]^{<\omega}|= |A|.$$
Therefore, the theorem now follows from Theorem \ref{Th:ConstructionOfModelsOfBAC}.
\Square 
\end{proof}        

Note that the Axiom of Choice in the metatheory is used in the above proof to show that that set of atoms is the same cardinality as the set of all finite sequences of atoms. If $\mathcal{M}= \langle M, \sqsubseteq^\mathcal{M} \rangle$ is countable, then the $\mathcal{N}= \langle N, \mathcal{S}^\mathcal{N}, \in^\mathcal{N}\rangle$ constructed in the proof of Theorem \ref{Th:EveryIABAIsASubsetRelation} is such that $\langle \mathcal{S}^\mathcal{N}, \subseteq^\mathcal{N} \rangle$ is the prime model of $\mathrm{Mer}$. The fact that the structure $\langle \mathcal{S}^\mathcal{N}, \subseteq^\mathcal{N} \rangle$ produced in the proof of Theorem \ref{Th:EveryIABAIsASubsetRelation} only contains finite joins of atoms of $\mathcal{M}$ prevents the sets of $\mathcal{N}$ from satisfying any form of the Axiom of Infinity and makes $\mathcal{N}$ look more like a model of second order arithmetic than class theory. The next result shows that every model of $\mathrm{BAS}$ can be realised as the sets of a model of $\mathrm{BAC}$.

\begin{Lemma1} \label{Th:ModelsOfMerCanBeExtendedToModelsOfIABA}
Let $\mathcal{M}= \langle M, \sqsubseteq^\mathcal{M} \rangle$ be a model of $\mathrm{Mer}$. Then there exists a model $\mathcal{N}= \langle N, \sqsubseteq^\mathcal{N} \rangle$ of $\mathrm{IABA}$ such that $\mathcal{M}$ is isomorphic to a proper ideal of $\mathcal{N}$ that contains all of the atoms of $\mathcal{N}$.   
\end{Lemma1}

\begin{proof}
Let $A= \{ x \in M \mid \mathcal{M} \models \mathbf{Atm}(x)\}$. For all $x \in M$, define $\mathrm{Ext}(x)$ to be the set $\{u \in A \mid u \sqsubseteq^\mathcal{M} x\}$. Let $X \subseteq \mathcal{P}(A)$ be the set 
$$X= \{ \mathrm{Ext}(x) \mid x \in M \}.$$
Let $f: M \longrightarrow X$ be defined by: for all $x \in M$, $f(x)= \mathrm{Ext}(x)$. Lemma \ref{Th:ElementsDeterminedByAtomsIABA} implies that $f$ witnesses the fact that $\mathcal{M} \cong \langle X, \subseteq \rangle$. Moreover, for all $x, y \in M$,
$$f\left((x \dot\lor y)^\mathcal{M}\right)= f(x) \cup f(y),\ f\left((x \dot\land y)^\mathcal{M}\right)= f(x) \cap f(y) \textrm{ and } f\left((x \dotminus y)^\mathcal{M}\right)= f(x) - f(y).$$
Let 
$$N= \{ x \subseteq A \mid (\exists y \in X)(x=y \lor x= A-y)\}.$$
Note that $\emptyset, A \in N$ and for all $x, y \in X$,
$$\begin{array}{cc}
x \cap (A-y)= x-y & x \cup (A-y)= A-(y-x)\\
(A-x) \cup (A-y)= A-(x\cap y) & (A-x) \cap (A-y)= A-(x \cup y)
\end{array}.$$
It follows that $N$ is closed under finite intersections, finite unions and complements. Therefore, $\langle N, \subseteq\rangle$ is a model of $\mathrm{IABA}$ with exactly the same atoms as $\langle X, \subseteq \rangle$. To see that $X$ is a proper ideal of $\langle N, \subseteq \rangle$, let $x, y \in X$. Note that if $A-y \subseteq x$, then $x \cup y= A$, which contradicts the fact that $A \notin X$. Therefore, $X$ is a proper ideal of $\langle N, \subseteq \rangle$.   
\Square   
\end{proof} 

\begin{Theorems1}
Let $\mathcal{M}= \langle M, \in^\mathcal{M} \rangle$ be a model of $\mathrm{BAS}$. Then there exists a model $\mathcal{N}= \langle N, \mathcal{S}^\mathcal{N}, \in^\mathcal{N} \rangle$ of $\mathrm{BAC}$ such that $\mathcal{M} \cong \langle \mathcal{S}^\mathcal{N}, \in^\mathcal{N}\rangle$. 
\end{Theorems1}

\begin{proof}
By Theorem \ref{Th:MainTheorem1}, $\langle M, \subseteq^\mathcal{M} \rangle$ is a model of $\mathrm{Mer}$ and $f= \{\langle x, \{x\} \rangle \mid x \in M\}$ witnesses the fact that $\langle M, \subseteq^\mathcal{M} \rangle$ has the same number of atoms as elements. Using Lemma \ref{Th:ModelsOfMerCanBeExtendedToModelsOfIABA}, let $\mathcal{N}= \langle N, \sqsubseteq^\mathcal{N} \rangle$ be a model of $\mathrm{IABA}$ with proper ideal $I \subseteq N$ containing all of the atoms of $\mathcal{N}$ and bijection $h: M \longrightarrow I$ witnessing the fact that $\langle M, \subseteq^\mathcal{M} \rangle \cong \langle I, \sqsubseteq^\mathcal{N} \rangle$. Note that $f^\prime= \{\langle h(x), h(\{x\}) \rangle \mid x \in M \}$ is a bijection witnessing the fact $\langle I, \sqsubseteq^\mathcal{N} \rangle$ has the same number of atoms as elements. Let $\in^* \subseteq N \times N \rangle$ be defined by: for all $x, y \in N$,
$$x \in^* y \textrm{ if and only if } \mathcal{N} \models f^\prime(x) \sqsubseteq y.$$
Therefore, by Theorem \ref{Th:ConstructionOfModelsOfBAC}, $\langle N, I, \in^* \rangle$ is model of $\mathrm{BAC}$ with $\subseteq^{\langle N, I, \in^* \rangle}= \sqsubseteq^\mathcal{N}$. Now, for all $x, y \in M$, 
$$\mathcal{M} \models x \in y  \textrm{ if and only if } \mathcal{M} \models f(x) \subseteq y \textrm{ if and only if } \mathcal{N} \models h(f(x)) \sqsubseteq h(y)$$ 
$$\textrm{if and only if } \mathcal{N} \models f^\prime(h(x)) \sqsubseteq h(y) \textrm{ if and only if } \langle N, I, \in^* \rangle \models h(x) \in h(y).$$
Therefore, $\langle I, \in^* \rangle \cong \mathcal{M}$.      
\Square
\end{proof}

We also get an analogue of Theorem \ref{Th:BASMinimal}:

\begin{Theorems1}
Let $\mathcal{M}= \langle M, \mathcal{S}^\mathcal{M}, \in^\mathcal{M} \rangle$ be a model of $\mathrm{CAS}$. Then the following are equivalent:
\begin{itemize}
\item[(I)] $\mathcal{M} \models \mathrm{BAC}$;
\item[(II)] $\langle M, \subseteq^\mathcal{M} \rangle \models \mathrm{IABA}$ and $\mathcal{S}^\mathcal{M}$ is a proper ideal of $\langle M, \subseteq^\mathcal{M} \rangle$. 
\end{itemize}
\end{Theorems1}

\begin{proof}
Note that $(I) \Rightarrow (II)$ follows from Theorem \ref{Th:SubsetRelationInBAC}. To see that $(II) \Rightarrow (I)$,  assume that $\langle M, \subseteq^\mathcal{M} \rangle$ is a model of $\mathrm{IABA}$ and $\mathcal{S}^\mathcal{M}$ is a proper ideal. The fact that $\langle M, \subseteq^\mathcal{M} \rangle$ is antisymmetric immediately implies that $\mathcal{M} \models \mathrm{CExt}$. Note that $\mathrm{Emp}$ and $\mathrm{Adj}$ imply that for all sets $x$, $\{x\}$ exists and is a set. Since every member of a class is a set ($\mathrm{Mem}$), this implies that the atoms of $\langle M, \subseteq^\mathcal{M} \rangle$ are exactly the singletons, which belong to $\mathcal{S}^\mathcal{M}$, and the function $f= \{\langle x, \{x\}\rangle \mid x \in \mathcal{S}^\mathcal{M} \}$ witnesses the fact that $\mathcal{S}^\mathcal{M}$ is the same size as the set of atoms of $\langle M, \subseteq^\mathcal{M} \rangle$. Define $\in^* \subseteq M \times M$ by: for all $x, y \in M$,
$$x \in^* y \textrm{ if and only if } \mathcal{M} \models f(x) \subseteq^\mathcal{M} y.$$
Therefore, by Theorem \ref{Th:ConstructionOfModelsOfBAC}, $\langle M, \mathcal{S}^\mathcal{M}, \in^* \rangle \models \mathrm{BAC}$. Now, for all $x, y \in M$,
$$x \in^* y \textrm{ if and only if } \{x\} \subseteq^\mathcal{M} y \textrm{ if and only if } x \in^\mathcal{M} y.$$
Therefore, $\mathcal{M} \models \mathrm{BAC}$.
\Square 
\end{proof}

Using the translation defined in Definition \ref{Df:TranslationOf2StratFormulae}, we obtain the following analogue of Theorem \ref{Th:TranslationOf2StratToSubsetAlgebra}.

\begin{Theorems1} \label{Th:ClassTranslationOf2StratToSubsetAlgebra}
Let $\phi(x_1, \ldots, x_n, y_1, \ldots, y_m)$ be an $\mathcal{L}$-formula with a stratification $\sigma: \mathbf{Var}(\phi) \longrightarrow 2$ such that for all $1 \leq i \leq n$, $\sigma(\textrm{`}x_i\textrm{'})= 0$ and for all $1 \leq i \leq m$, $\sigma(\textrm{`}y_i\textrm{'})= 1$. Let $\mathcal{M}= \langle M, \mathcal{S}^\mathcal{M}, \in^\mathcal{M} \rangle$ be a model of $\mathrm{BAC}$. For all $a_1, \ldots, a_n, b_1, \ldots, b_m \in M$,
$$\langle M, \in^\mathcal{M} \rangle \models \phi(a_1, \ldots, a_n, b_1, \ldots, b_m) \textrm{ if and only if}$$
$$\langle M, \subseteq^\mathcal{M} \rangle \models \phi^\tau(\{a_1\}, \ldots, \{a_n\}, b_1, \ldots, b_m),$$
where $\tau$ is the translation of $2$-stratified $\mathcal{L}$-formulae to $\mathcal{L}_\mathrm{po}$-formulae described in Definition \ref{Df:TranslationOf2StratFormulae}.\Square  
\end{Theorems1}

Combined with Theorem \ref{Th:IABAComplete} this yields:

\begin{Coroll1} \label{Th:BACDecides2StratSentences}
If $\phi$ is a $2$-stratified $\mathcal{L}$-sentence, then 
$$\mathrm{BAC} \vdash \phi \textrm{ or } \mathrm{BAC} \vdash \neg \phi.$$
\Square 
\end{Coroll1} 

\begin{Coroll1} \label{Th:2StratCosequencesOfBACRecursive}
Let $T \supseteq \mathrm{BAC}$ be consistent. Then
$$\{\ulcorner \phi \urcorner \mid (\phi \textrm{ is a } 2 \textrm{-stratified } \mathcal{L}\textrm{-sentence}) \land (T \vdash \phi) \} \textrm{ is recursive.}$$
\Square
\end{Coroll1}

Note that that Theorem \ref{Th:ClassTranslationOf2StratToSubsetAlgebra}, and Corollaries \ref{Th:BACDecides2StratSentences} and \ref{Th:2StratCosequencesOfBACRecursive} only refer to $2$-stratified $\mathcal{L}$-sentences. The proof of the following is based on \cite[Remark 2.4]{dm13}:

\begin{Theorems1}
The $2$-stratified $\mathcal{L}_\mathrm{cl}$-sentence
$$\phi: \forall x \exists y (\forall z(z \in y \iff z \notin x) \land (\mathcal{S}(x) \lor \mathcal{S}(y)))$$
is independent of $\mathrm{BAC}$.  
\end{Theorems1}

\begin{proof}
Let 
$$M= \{X \subseteq \omega \mid |X| < \omega \lor |\omega \backslash X|< \omega \},$$ 
$$N= \mathcal{P}(\omega) \textrm{ and } I= \{ X \subseteq \omega \mid |X| < \omega \}.$$
Note that $\langle M, \subseteq \rangle$ and $\langle N, \subseteq \rangle$ both satisfy $\mathrm{IABA}$. Moreover,
$$I \supseteq \{\{n\} \mid n \in \omega\}= \{x \in M \mid \langle M, \subseteq \rangle \models \mathbf{Atm}(x) \}= \{x \in N \mid \langle N, \subseteq \rangle \models \mathbf{Atm}(x)\},$$
$I$ is a proper ideal of $\langle M, \subseteq \rangle$ and $\langle N, \subseteq \rangle$, and $|I|= \omega$. Applying Theorem \ref{Th:ConstructionOfModelsOfBAC}, we obtain $\in^M \subseteq M \times M$ and $\in^N \subseteq N \times N$ such that 
$$\langle M, I, \in^M \rangle \models \mathrm{BAC} \textrm{ and } \langle N, I, \in^N \rangle \models \mathrm{BAC}.$$
Now,
$$\langle M, I, \in^M \rangle \models \phi \textrm{ and } \langle N, I, \in^N \rangle \models \neg \phi.$$
\Square   
\end{proof}

We now turn identifying the complete theory that holds in every structure \mbox{$\langle M, \mathcal{S}^\mathcal{M}, \subseteq^\mathcal{M} \rangle$} when $\mathcal{M}= \langle M, \mathcal{S}^\mathcal{M}, \in^\mathcal{M} \rangle$ is a model of a sufficiently strong subsystem of $\mathrm{NBG}$. Let $\mathcal{L}_{\mathrm{Ideal}}$ be the language obtained from $\mathcal{L}_{\mathrm{po}}$ by adding a unary predicate $\mathcal{I}$. Therefore an $\mathcal{L}_{\mathrm{Ideal}}$-structure is a triple $\mathcal{M}= \langle M, \mathcal{I}^\mathcal{M}, \sqsubseteq^\mathcal{M} \rangle$ with $\mathcal{I}^\mathcal{M} \subseteq M$. The following extension of $\mathrm{IABA}$ appears as $T_1$ in \cite[Example 2]{dm13}:
\begin{itemize}
\item $\mathrm{IABA}_{\mathrm{Ideal}}$ is the $\mathcal{L}_{\mathrm{Ideal}}$-theory extending $\mathrm{IABA}$ with axioms:
\begin{itemize}
\item[](Ideal Axioms) Axioms saying that $\mathcal{I}(x)$ defines a proper ideal. And, for all $n \in \mathbb{N}$, the axiom that asserts that for all $x$, if $x$ has at most $n$ atoms below it, then $\mathcal{I}(x)$. I.e. Axioms saying that $\mathcal{I}(x)$ defines a proper ideal that contains the proper ideal of all elements that are above only finitely many atoms.
\item[](Main Axiom)
$$\forall x(\neg \mathcal{I}(x) \Rightarrow \exists y(y \sqsubseteq x \land \neg \mathcal{I}(y) \land \neg \mathcal{I}(x \dot\land \dotminus y) ))$$  
\end{itemize}
\end{itemize}
   
Derakhshan and MacIntyre \cite[Theorem 2.3]{dm13} show that the theory $\mathrm{IABA}_{\mathrm{Ideal}}$ admits quantifier elimination in an expansion of the language $\mathcal{L}_{\mathrm{Ideal}}$ that includes, for all $n \geq 1$, a unary predicate $C_n(x)$ whose interpretation is that there are at least $n$ atoms below $x$. 

\begin{Theorems1}
(Derakhshan-MacIntyre) The theory $\mathrm{IABA}_{\mathrm{Ideal}}$ is complete and decidable. \Square 
\end{Theorems1}

We now extend the theory $\mathrm{BAC}$ to a theory $\mathrm{BAC}^+$ with the property that if $\mathcal{M}= \langle M, \mathcal{S}^\mathcal{M}, \in^\mathcal{M} \rangle$ is a model of $\mathrm{BAC}^+$, then the structure $\langle M, \mathcal{S}^\mathcal{M}, \subseteq^\mathcal{M} \rangle$ is an $\mathcal{L}_\mathrm{Ideal}$-structure that satisfies $\mathrm{IABA}_{\mathrm{Ideal}}$.
\begin{itemize}
\item $\mathrm{BAC}^+$ is the $\mathcal{L}_\mathrm{cl}$-theory extending $\mathrm{BAC}$ with the axiom:
$$\mathrm{(Sep)} \ \ \forall X \left( \begin{array}{c}
\neg \exists x(X=x) \Rightarrow\\
\exists Y \exists Z\left(\begin{array}{c}
\forall w(w \in Y \Rightarrow w \in X) \land\\
\forall v(v \in Z \iff (v \in X \land v \notin Y)) \land\\
\neg \exists y(y=Y) \land \neg z(z=Z))
\end{array}\right)   
\end{array}\right)$$
\end{itemize}
Note that, using set-theoretic operations available in $\mathrm{BAC}$, the axiom $\mathrm{Sep}$ can be written as $\forall X(\neg \mathcal{S}(X) \Rightarrow \exists Y(Y \subseteq X \land \neg \mathcal{S}(Y) \land \neg \mathcal{S}(X \backslash Y)))$. This observation, combined with Theorems \ref{Th:SubsetRelationInBAC} and \ref{Th:ConstructionOfModelsOfBAC}, yields:

\begin{Theorems1} \label{Th:TheoryOfSubsetOrderWithSInBACPlus}
Let $\mathcal{M}= \langle M, \mathcal{S}^\mathcal{M}, \in^\mathcal{M} \rangle$ be a model of $\mathrm{BAC}^+$. Then $\langle M, \mathcal{S}^\mathcal{M}, \subseteq^\mathcal{M} \rangle$ is a model of $\mathrm{IABA}_{\mathrm{Ideal}}$ and $\mathcal{S}^\mathcal{M}$ is the same size as the set of atoms of $\langle M, \subseteq^\mathcal{M} \rangle$. \Square
\end{Theorems1}

\begin{Theorems1}
Let $\mathcal{M}= \langle M, \mathcal{I}^\mathcal{M}, \sqsubseteq^\mathcal{M} \rangle$ be a model of $\mathrm{IABA}_{\mathrm{Ideal}}$ and let $f: \mathcal{I}^\mathcal{M} \longrightarrow A$ be a bijection, where $A= \{ x \in M \mid \mathcal{M} \models \mathbf{Atm}(x)\}$. Define $\in^* \subseteq M \times M$ by: for all $x, y \in M$,
$$x \in^* y \textrm{ if and only if } \mathcal{M} \models f(x) \sqsubseteq y.$$
Then $\langle M, \mathcal{I}^\mathcal{M}, \in^* \rangle \models \mathrm{BAC}^+$ and $\subseteq^{\langle M, \mathcal{I}^\mathcal{M}, \in^* \rangle}= \sqsubseteq^\mathcal{M}$. \Square
\end{Theorems1}

\begin{Theorems1}
Let $\mathcal{M}= \langle M, \mathcal{S}^\mathcal{M}, \in^\mathcal{M} \rangle$ be a model of $\mathrm{CAS}$. Then the following are equivalent:
\begin{itemize}
\item[(I)] $\mathcal{M} \models \mathrm{BAC}^+$;
\item[(II)] $\langle M, \mathcal{S}^\mathcal{M}, \subseteq^\mathcal{M} \rangle \models \mathrm{IABA}_{\mathrm{Ideal}}$.
\end{itemize}
\Square
\end{Theorems1}

Modifying Example \ref{Ex:ModelOfMerThatCannotBeSubset}, we exhibit a model $\mathrm{IABA}_{\mathrm{Ideal}}$ that is not isomorphic to a structure $\langle M, \mathcal{S}^\mathcal{M}, \subseteq^\mathcal{M} \rangle$ arising from a model $\mathcal{M}= \langle M, \mathcal{S}^\mathcal{M}, \in^\mathcal{M} \rangle$ of $\mathrm{BAC}^+$.

\begin{Examp1}
Let $M= \mathcal{P}(\omega+\omega)$ and let $I=\{X \in M \mid (\exists \beta \in \omega+\omega)(\forall \alpha \in X)(\alpha \in \beta)\}$. Now, $\langle M, I, \subseteq \rangle$ is a model of $\mathrm{IABA}_{\mathrm{Ideal}}$. But,
$$|\{x \in M \mid \langle M, I, \subseteq \rangle \models \mathbf{Atm}(x) \}|= |\{\{\alpha\} \mid \alpha \in \omega+\omega \}|= \aleph_0$$
and $|I|= 2^{\aleph_0}$. So, by Theorem \ref{Th:TheoryOfSubsetOrderWithSInBACPlus}, there does not exist $\in^* \subseteq M \times M$ such that $\langle M, I, \in^* \rangle$ is a model of $\mathrm{BAC}^+$.     
\end{Examp1}

The realisations of models of $\mathrm{Mer}$ as ideals of models of $\mathrm{IABA}$ obtained in Lemma \ref{Th:ModelsOfMerCanBeExtendedToModelsOfIABA} does not yield models of $\mathrm{IABA}_{\mathrm{Ideal}}$.

\begin{Quest1}
Can every model $\langle I, \in^* \rangle$ of $\mathrm{BAS}$ be expanded to a model $\langle M, I, \in^* \rangle$ of $\mathrm{BAC}^+$?
\end{Quest1} 

We next turn to showing that the theory $\mathrm{BAC}^+$ is a subsystem of $\mathrm{NBG}$. In the theory $\mathrm{NBG}$, we use $\mathrm{Ord}$ to denote the class of sets that are von Neumann ordinals. 

\begin{Lemma1} \label{Th:PropertiesOfProperClassesNBG}
The theory $\mathrm{NBG}$ proves that for all $X$, the following are equivalent:
\begin{itemize}
\item[(i)] $\neg \mathcal{S}(X)$;
\item[(ii)] there exists $f: X \longrightarrow \mathrm{Ord}$ such that the range of $f$ is unbounded in $\mathrm{Ord}$;
\item[(iii)] there exists a surjection $f: X \longrightarrow \mathrm{Ord}$. 
\end{itemize}
\end{Lemma1}

\begin{proof}
$(i) \Rightarrow (ii)$: Assume that $X$ is not a set. The restriction of the usual set-theoretic rank function to $X$ defines a a function $f:X \longrightarrow \mathrm{Ord}$. If the range of this function is bounded in $\mathrm{Ord}$, then $X$ is a set.\\
$(ii) \Rightarrow (iii)$: Let $f: X \longrightarrow \mathrm{Ord}$ be such that $Y= \mathrm{rng}(f)$ is unbounded in $\mathrm{Ord}$. Let $g: Y \longrightarrow \mathrm{Ord}$ be the inverse of the order preserving function that enumerates the elements of $Y$. If $g$ is not surjective, then the range of $g$ is bounded in $\mathrm{Ord}$. But, if the range of $g$ is bounded in $\mathrm{Ord}$, then Replacement in $\mathrm{NBG}$ applied to $g^{-1}$ implies that $Y$ is a set, which is a contradiction. Therefore, $g$ is surjective and $g \circ f$ is the desired function.\\
$(iii) \Rightarrow (i)$: Let $g: X \longrightarrow \mathrm{Ord}$ be a surjection. If $X$ is a set, then Replacement in $\mathrm{NBG}$ implies that $\mathrm{Ord}$ is set, which is a contradiction.  
\Square
\end{proof}

\begin{Theorems1}
The theory $\mathrm{BAC}^+$ is a subtheory of $\mathrm{NBG}$.  
\end{Theorems1}

\begin{proof}
It is clear that $\mathrm{BAC}$ is a subtheory of $\mathrm{NBG}$, so we only need to verify that $\mathrm{NBG}$ proves the $\mathrm{Sep}$ axiom. Let $X$ be such that $\neg \mathcal{S}(X)$. Using Lemma \ref{Th:PropertiesOfProperClassesNBG}, let $f: X \longrightarrow \mathrm{Ord}$ be surjective. Let $Y= \{ x \in X \mid f(x) \textrm{ is a limit ordinal}\}$. Now, the ranges of both $f \upharpoonright Y$ and $f \upharpoonright (X \backslash Y)$ are unbounded in $\mathrm{Ord}$. Therefore, by Lemma \ref{Th:PropertiesOfProperClassesNBG}, $\neg \mathcal{S}(Y)$ and $\neg \mathcal{S}(X \backslash Y)$. This shows that $\mathrm{NBG}$ proves $\mathrm{Sep}$. 
\Square
\end{proof}

\begin{Coroll1}
If $\mathcal{M}= \langle M, \mathcal{S}^\mathcal{M}, \in^\mathcal{M} \rangle$ is a model of $\mathrm{NBG}$, then $\langle M, \mathcal{S}^\mathcal{M}, \subseteq^\mathcal{M} \rangle \models \mathrm{IABA}_{\mathrm{Ideal}}$. \Square 
\end{Coroll1}

Theorem \ref{Th:TheoryOfSubsetOrderWithSInBACPlus} allows us to extend the analysis of section \ref{Sec:2StratifiedSentences} to show that the theory $\mathrm{BAC}^+$ decides every $2$-stratified $\mathcal{L}_\mathrm{cl}$-sentence. 

\begin{Definitions1}
Let $\phi$ be an $\mathcal{L}_\mathrm{cl}$-formula with a stratification $\sigma: \mathbf{Var}(\phi) \longrightarrow 2$. Define a translation, $\chi$, of $\phi$ into a $\mathcal{L}_{\mathrm{Ideal}}$-formula $\phi^\chi$ by recursion using the rules described in Definition \ref{Df:TranslationOf2StratFormulae} supplemented with the rule: if $\phi$ is $\mathcal{S}(u)$, then $\phi^\chi$ is $\mathcal{I}(u)$.
\end{Definitions1}

The translation $\chi$ yields the following analogue of Theorem \ref{Th:TranslationOf2StratToSubsetAlgebra}:

\begin{Theorems1}
Let $\phi(x_1, \ldots, x_n, y_1, \ldots, y_m)$ be an $\mathcal{L}_\mathrm{cl}$-formula with a stratification $\sigma: \mathbf{Var}(\phi) \longrightarrow 2$ such that for all $1 \leq i \leq n$, $\sigma(\textrm{`}x_i\textrm{'})= 0$, and for all $1 \leq i \leq m$, $\sigma(\textrm{`}y_i\textrm{'})= 1$. Let $\mathcal{M}= \langle M, \mathcal{S}^\mathcal{M}, \in^\mathcal{M} \rangle$ be a model of $\mathrm{BAC}^+$. For all $a_1, \ldots, a_n, b_1, \ldots, b_m \in M$,
$$\mathcal{M} \models \phi(a_1, \ldots, a_n, b_1, \ldots, b_m) \textrm{ if and only if}$$
$$\langle M, \mathcal{S}^\mathcal{M}, \subseteq^\mathcal{M} \rangle \models \phi^\chi(\{a_1\}, \ldots, \{a_n\}, b_1, \ldots, b_m).$$
\Square
\end{Theorems1}

This shows that $\mathrm{BAC}^+$ decides every $2$-stratified $\mathcal{L}_\mathrm{cl}$-sentence.

\begin{Coroll1}
If $\phi$ is a $2$-stratified $\mathcal{L}_\mathrm{cl}$-sentence, then 
$$\mathrm{BAC}^+ \vdash \phi \textrm{ or } \mathrm{BAC}^+ \vdash \neg \phi.$$
\Square
\end{Coroll1}

\begin{Coroll1}
Let $T \supseteq \mathrm{BAC}^+$ be consistent. Then 
$$\{ \ulcorner \phi \urcorner \mid (\phi \textrm{ is a } 2\textrm{-stratified } \mathcal{L}_\mathrm{cl}\textrm{-sentence}) \land (T \vdash \phi) \} \textrm{ is recursive.}$$
\Square
\end{Coroll1}

\bibliographystyle{alpha}
\bibliography{.}

\begin{thebibliography}{9}

\bibitem[Che]{cheXX} Cheng, Y. ``Current research on G\"{o}del's Incompleteness Theorems". {\it Bulletin of Symbolic Logic}. Vol. 27. 2021. pp 113-167.

\bibitem[DM]{dm13} Derakhshan, J. and MacIntyre, A. ``Enrichments of Boolean algebras: A uniform treatment of some classical and some novel examples". Available online at {\tt arXiv.org} ({\tt arXiv:1310.3527}).  

\bibitem[Er\v{s}]{ers64} Er\v{s}ov, J. L. ``Decidability of the elementary theory of relatively complemented lattices and of the theory of filters". {\it Algebra i Logika}. Vol. 3. No. 3. 1964. pp 17-38.

\bibitem[For]{for92} Forster, T. E. {\it Set theory with a Universal set}. Oxford Logic Guides No. 31. Oxford University Press. 1992.

\bibitem[Fri]{fri03} Friedman, H. ``Three quantifier sentences". {\it Fundamenta Mathematicae}. Vol. 177. 2003. pp 213-240.

\bibitem[Gog78]{gog78} Gogol, D. ``The $\forall_n \exists$-Completeness of Zermelo-Fraenkel Set Theory". {\it Zeitschrift f\"{u}r mathematische Logik und Grundlagen der Mathematik}. Vol. 24. 1978. pp 289-290.

\bibitem[Gog79]{gog79} Gogol, D. ``Sentences with three quantifiers are decidable in set theory". {\it Fundamenta Mathematicae}. CII. 1979. pp 1-8.

\bibitem[Gri]{gri73} Gri\u{s}hin, V. N. {\it The method of stratification in set theory}. Academy of Sciences of the USSR. Moscow, 1973. 

\bibitem[HK16]{hk16} Hamkins, J. D. and Kikuchi, M. ``Set-theoretic mereology". {\it Logic and Logical Philosophy Special Issue, special issue ``Mereology and beyond, part II"}. Vol. 25. No. 3. 2016. pp 285-308.

\bibitem[HK17]{hkxx} Hamkins, J. D. and Kikuchi, M. ``The inclusion relations of the countable models of set theory are all isomorphic". Available online at {\tt arXiv.org} ({\tt arXiv:1704.04480}).

\bibitem[Men]{men97} Mendelson, E. {\it Introduction to mathematical logic}. Fourth edition. Chapman and Hall, London, 1997.  

\bibitem[Qui]{qui37} Quine, W. v. O. ``New Foundations for mathematical logic". {\it American Mathematical Monthly}. Vol. 44. 1937. pp 70-80.

\bibitem[Tar]{tar49} Tarski, A. ``Arithmetical classes and types of Boolean algebras". {\it Bulletin of the American Mathematical Society}. Vol. 55. 1949. p 63.

\bibitem[Vis]{vis10} Visser, A. ``What is the right notion of sequentiality?". {\it Logic preprint series}. Vol. 288. 2010. pp 1-24. 

\end{thebibliography}

\end{document}